\documentclass[epj,referee]{svjour}
\usepackage{graphics}
\usepackage{graphicx}
\usepackage{epsfig}
\usepackage{amssymb}

\begin{document}
\title{Evolutionary Prisoner's Dilemma on heterogeneous Newman-Watts small-world network}
\author{F. Fu\thanks{\email{fufeng@pku.edu.cn}} \and L.-H. Liu \and
L. Wang\thanks{\email{longwang@pku.edu.cn}}
}                     
%
%
\institute{Center for Systems and Control, College of Engineering,
Peking University, Beijing 100871, P.~R. China}

\date{Received: date / Revised version: date}
%
\abstract{In this paper, we focus on the heterogeneity of social
networks and its role to the emergence of prevailing cooperation
and sustaining cooperators. The social networks are representative
of the interaction relationships between players and their
encounters in each round of games. We study an evolutionary
Prisoner's Dilemma game on a variant of Newman-Watts small-world
network, whose heterogeneity can be tuned by a parameter. It is
found that optimal cooperation level exists at some intermediate
topological heterogeneity for different temptations to defect.
That is, frequency of cooperators peaks at some specific values of
degree heterogeneity --- neither the most heterogeneous case nor
the most homogeneous one would favor the cooperators. Besides, the
average degree of networks and the adopted update rule also affect
the cooperation level.
} 
\maketitle
\section{Introduction}
Evolutionary game theory has been well developed as an
interdisciplinary science by researchers from biology, economics,
social science, computer science for several decades. In past few
years, it also gained the interests of physicists to study some
phenomena and intriguing mechanisms in well-mixed population by
using mean-field theory of statistical physics. In classic game
theory, the players are assumed to be completely rational and try
to maximize their utilities according to opponents' strategies.
Under these assumptions, in the Prisoner's Dilemma game (PD game),
two players simultaneously decide whether to cooperate (C) or to
defect (D). They both receive $R$ upon mutual cooperation and $P$
upon mutual defection. A defector exploiting a C player gets  $T$,
and the exploited cooperator receives $S$, such that $T>R>P>S$ and
$2R>T+S$. As a result, it is best to defect regardless of the
co-player's decision. Thus, defection is the evolutionarily stable
strategy (ESS), even though all individuals would be better off if
they cooperated. Thereby this creates the social dilemma, because
when everybody defects, the mean population payoff is lower than
that when everybody cooperates. However, cooperation is ubiquitous
in natural systems from cellular organisms to mammals. In the past
two decades, some extensions on PD game have been considered to
elucidate the underlying mechanisms boosting cooperation behaviors
by which this dilemma could be resolved. For instance, depart from
the well-mixed population scenario, Nowak and May considered PD
game on spatially structured population~\cite{Nowak92}. All
individuals are constrained to play only with their immediate
neighbors. They found that spatial structure enhances the ability
of cooperators to resist invasion by defectors.

The successful development of network science provides a
systematic framework for studying the dynamics processes taking
place on complex networks \cite{Albert2002,Newman2003}. For most
networks, including the World Wide Web, the Internet, and the
metabolic networks, they have small-world property and are
demonstrated to be scale-free~\cite{Watts98,Barabasi99}. The
network theory is a natural and convenient tool to describe the
population structure on which the evolution of cooperation is
studied. Each vertex represents an individual and the edges denote
links between players in terms of game dynamical interaction.
Trivially, one could conclude that both well-mixed population and
spatially-structured population are represented by regular graphs.
It is, however, unrealistic to assume the real world is as
homogeneous as regular graphs. Furthermore, as aforementioned,
many real-world networks of interactions are heterogeneous,
namely, different individuals have different numbers of average
neighbors with whom they interact. In~\cite{Abra2001}, Abramson
and Kuperman studied a simple model of an evolutionary version of
the Prisoner's Dilemma game played in small-world networks. They
found that defectors prevail at some intermediate rewiring
probability. By introducing the volunteering participation,
Szab\'{o} {\it et al.} investigated the spatio-temporal diagrams
of the evolutionary process within the context of rewired lattices
which have small-world property~\cite{Szabo2004}. Other paradigms
of game with spatially-structured population or on graphs, such as
snowdrift game, rock-scissors-paper game, public goods game, etc,
have been studied~\cite{Toma2006,Szolnoki2004,Szabo2002prl}.

Recently, Santos and Pacheco discovered that in scale-free
networks, the cooperators become dominant for entire range of
parameters in evolutionary PD game and snowdrift
game~\cite{Santos2005}. They also demonstrated that the
enhancement of cooperation would be inhibited whenever the
correlations between individuals are decreased or removed. Being
similar to Hamilton's rule, Ohtsuki {\it et al} described a simple
rule for evolution of cooperation on graphs and social networks
--- the benefit-to-cost ratio exceeds the average number of
neighbors~\cite{Ohtsuki2006}. In~\cite{Santos2006a,Santos2006b},
Santos suggested that the heterogeneity of the population
structure, especially scale-free networks, provides a new
mechanism for cooperators to survive. Moreover, cooperators are
capable of exploiting the heterogeneity of the network of
interactions, namely, finally occupy the hubs of the networks.
Besides, Vukov and  Szabo studied the evolutionary PD game on
hierarchical lattices and revealed that the highest stationary
frequency of cooperators occurs in some intermediate
layers~\cite{Vukov2005}. Some researchers have studied minority
game~\cite{YCZ97} , with small-world interactions~\cite{Kirley06},
with emergent scale-free leadership structure~\cite{Angehel04},
etc. Howbeit, these works are mostly studied on crystalized
(static) networks, i.e., the topology of network is not affected
by the dynamics on it. Therefore, adapting the network topology
dynamically in response to the dynamic state of the nodes should
be taken into account. Evolutionary games with ``adaptive''
network structure have been investigated by some
researchers~\cite{Zim05,Ren2006}.

In this paper, we consider evolutionary PD game on a variant of
Newman-Watts (NW) small-world network. As is known, the NW
small-world network is moderately homogeneous.  In order to
investigate the heterogeneity's role in the emergence of
cooperation for different temptations to defect, some artificial
heterogeneity in network degree sequences is introduced. Namely,
some nodes are randomly chosen as hubs and at least one endpoint
of each of the added $m$ shortcuts is linked to these hubs. The
remainder of the paper is organized as follows. Section II
discusses the heterogeneous Newman-Watts (HNW) model and the
evolutionary PD game, including payoff matrix and the microscopic
updating rule. And then Sec. III gives out the simulation results
for different parameters and makes some explanations to the
results. Conclusions are made in Sec. IV.

\section{The model}

\subsection{The Heterogeneous Newman-Watts (HNW)  model}
\begin{figure} \centering
\includegraphics[width=5cm]{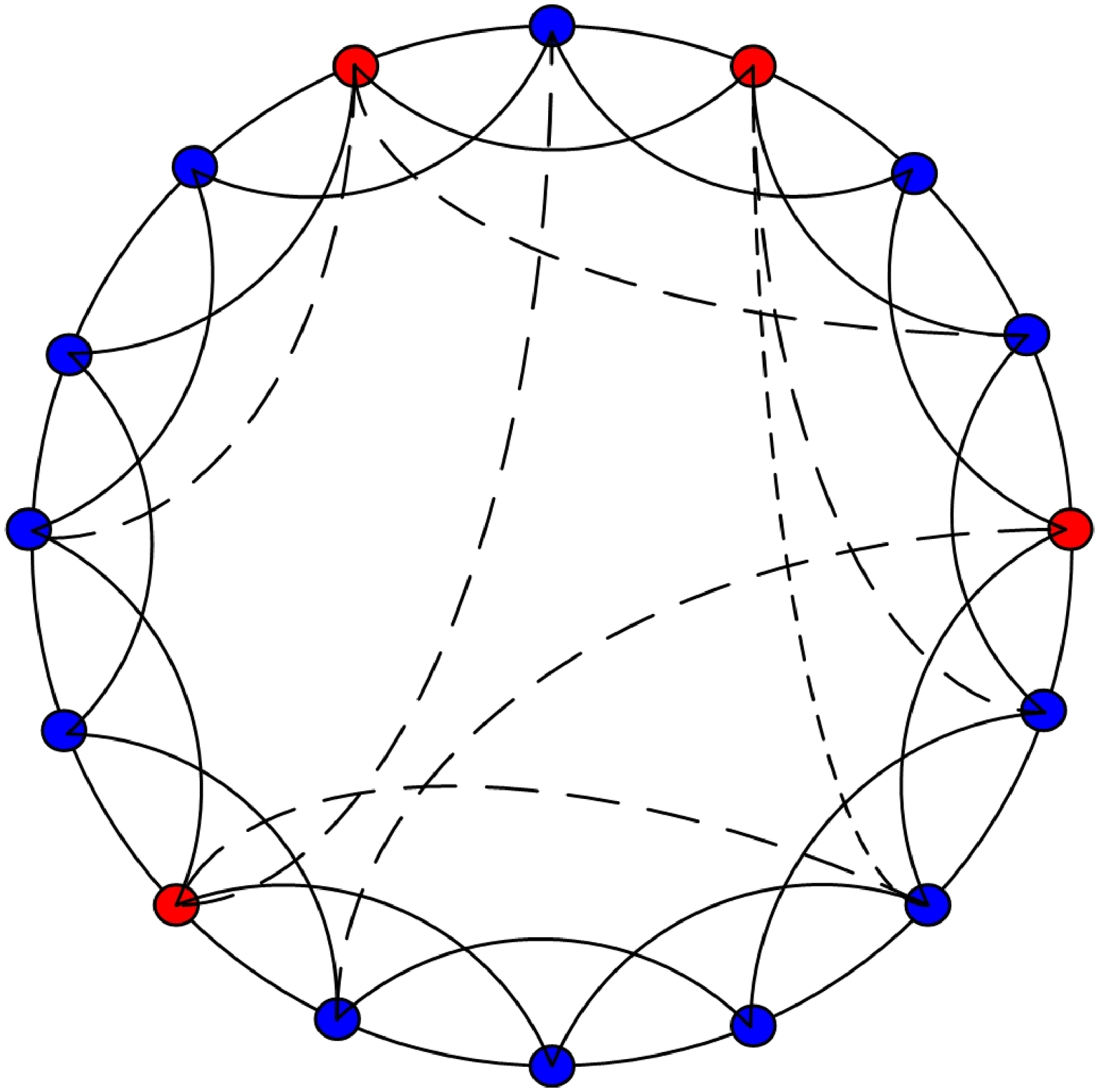}
\caption{\label{fig1}}
\end{figure}

\begin{figure} \centering
\includegraphics[width=7.5cm]{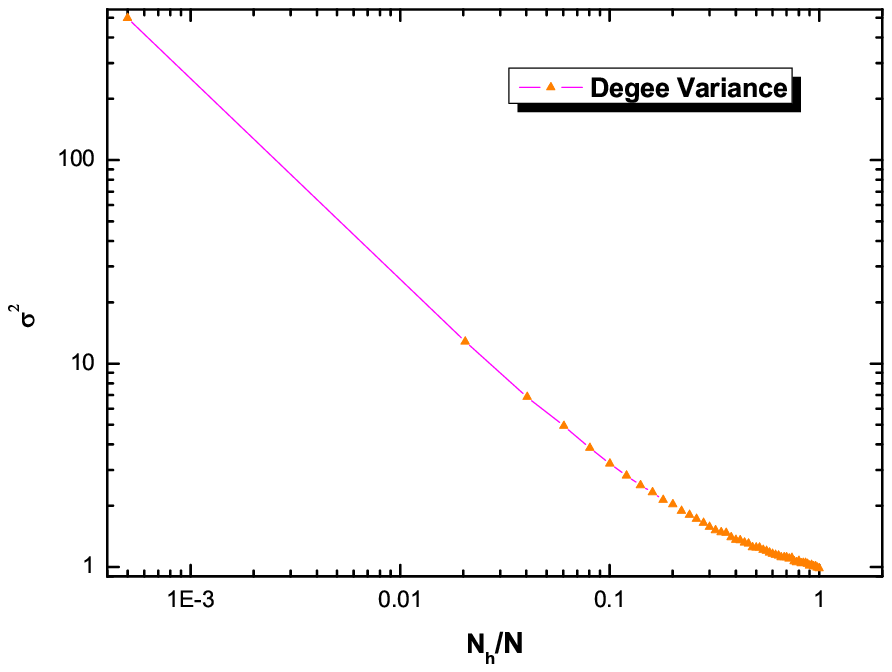}
\caption{\label{fig2}}
\end{figure}

We consider one-dimensional lattice of $N$ vertices with periodic
boundary conditions, i.e., on a ring, and join each vertex to its
neighbors $\kappa$ or fewer lattice spacing away. Instead of
rewiring a fraction of the edges in the regular lattice as
proposed by Watts and Strogatz~\cite{Watts98}, some ``shortcuts''
connecting randomly chosen vertex pairs are added to the
low-dimensional lattice~\cite{NewmanWatts}. To introduce certain
type of heterogeneity to this NW model~\cite{Nishikawa03}, we
choose $N_h$ vertices at random from all $N$ nodes with equal
probabilities and treat them as hubs. Then we add each of $m$
shortcuts by connecting one node at random from all $N$ nodes to
another node randomly chosen from the $N_h$ centered nodes
(duplicate connections and self-links are voided), see
Fig.~\ref{fig1} as an illustration. In such heterogeneous
Newman-Watts (HNW) model, nodes can be naturally divided into two
groups: hubs, which tend to have higher connectivity, and the
others that have lower connectivity. Accordingly, the
heterogeneity is controlled by the parameter $N_h$ : small $N_h$
leads to higher degree of the hubs, which in turn results in
increased heterogeneity. Trivially, one could see that when
$N_h=N$, shortcuts are just simply added uniformly at random,
effectively reducing the network to the homogeneous Newman-Watts
model. If $N_h=1$, all shortcuts are linked to a single center,
making the network an extremely heterogeneous one. The degree of
heterogeneity of the network is computed as
$h=N^{-1}\sum_kk^2N(k)-\langle k \rangle$ (the variance of the
network degree sequcence), where $N(k)$ gives the number of
vertices with $k$ edges, $\langle k \rangle$ denotes the average
degree. From Fig.~\ref{fig2}, it shows that variance
$\sigma^2(=h)$ monotonically decreases as the fraction number of
hubs $N_h/N$ increases. Thus we simply use the fraction $N_h/N$ to
indicate the degree of heterogeneity.

\subsection{Evolutionary Prisoner's Dilemma (PD) game}
Since the pioneering work on iterated games by
Axelrod~\cite{Axelrod81}, the evolutionary Prisoner's Dilemma (PD)
game has been a general metaphor for studying the cooperative
behavior. In the evolutionary PD game, the individuals are pure
strategists, following two simple strategies: cooperate (C) and
defect (D). The spatial distribution of strategies is described by
a two-dimensional unit vector for each player $x$, namely,
\begin{equation}
s=\left(\begin{array}{c}
  1 \\
  0 \\
\end{array}\right )\,\, \mbox{and}\,\,
\left(\begin{array}{c}
  0 \\
  1 \\
\end{array}\right )
\end{equation}
for cooperators and defectors, respectively. Each individual plays
the PD game with its ``neighbors'' defined by their who-meets-whom
relationships and the incomes are accumulated. The total income of
the player at the site $x$ can be expressed as
\begin{equation}
P_x=\sum_{y\in \Omega_x}s_x^TMs_y
\end{equation}
where the $\Omega_x$ denotes the neighboring sites of $x$, and the
sum runs over neighbor set $\Omega_x$ of the site $x$. Without
loss of the generic feature of PD game, the payoff matrix has a
rescaled form, suggested by Nowak and May~\cite{Nowak93}
\begin{equation}
M=\left (\begin{array}{cc}
  1 \,& \,0 \\
  b \,& \,0 \\
\end{array}\right )
\end{equation}
where $1<b<2$. In evolutionary games the players are allowed to
adopt the strategy of one of their more successful neighbors. Then
the individual $x$ randomly selects a neighbor $y$ for possibly
updating its strategy. Since the success is measured by the total
payoff, whenever $P_y>P_x$ the site $x$ will adopt $y$'s strategy
with probability given by~\cite{Santos2005}
\begin{equation}
\label{transp} W_{s_x \leftarrow s_y} = \frac{P_y-P_x}{bk_>}
\end{equation}
where $k_>$ is the largest between site $x$'s degree $k_x$ and
$y$'s $k_y$. This microscopic updating rule is some kind of
imitation process that is similar to Win-Stay-Lose-Shift rules in
spirit. Unlike the case of bounded rationality, such proportional
imitation (the take-over probability is proportional to the payoff
difference of the two sites) does not allow for an inferior
strategy to replace a more successful one.

\section{Simulation results and discussions}
\begin{figure} \centering
\includegraphics[width=7.5cm]{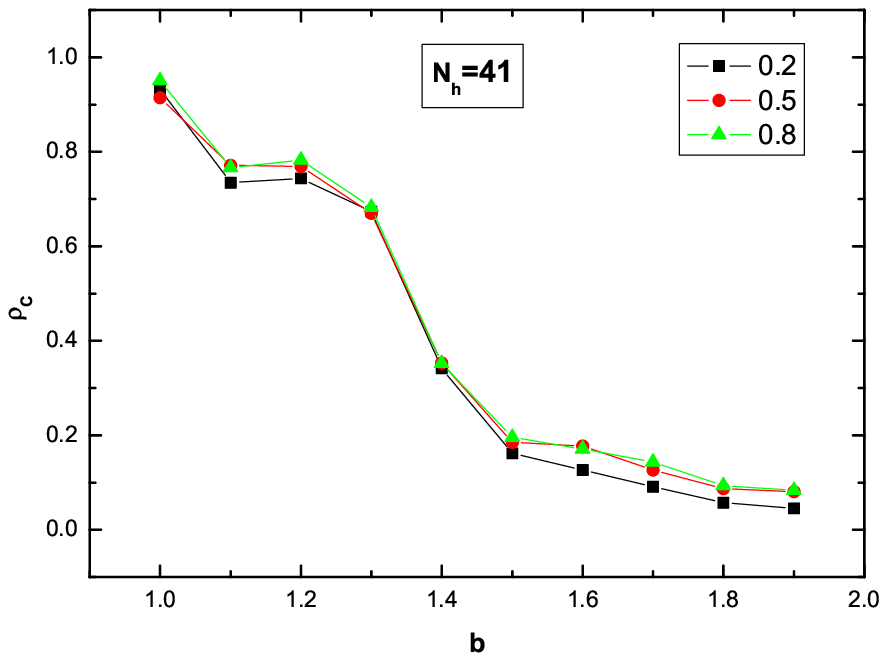}
\caption{\label{fig3}}
\end{figure}

Evolutionary PD game is performed analogically to replicator
dynamics: in each generation, all directly connected pairs of
individuals $x$ and $y$ engage in a single round game and their
accumulated payoffs are denoted by $P_x$ and $P_y$ respectively.
The synchronous updating rule is adapted here. Whenever a site $x$
is updated, a neighbor $y$ is chosen at random from all $k_x$
neighbors. The chosen neighbor takes over site $x$ with
probability $W_{s_x \leftarrow s_y}$ as Eq.~(\ref{transp}),
provided that $P_y>P_x$. Simulations were carried out for a
population of $N=2001$ players occupying the vertices of the HNW
network. Initially, an equal percentage of cooperators and
defectors was randomly distributed among the population. We
confirm that different initial conditions do not qualitatively
influence the equilibrium results. In Fig.~\ref{fig3}, with fixed
chosen hub numbers $N_h=41$, the equilibrium frequencies of
cooperators with respect to different $b$ are almost the same when
started from different initial conditions. Equilibrium frequencies
of cooperators were obtained by average over 2000 generations
after a transient time of 10000 generations. The evolution of the
frequency of cooperators as a function of $b$ and $N_h/N$ has been
computed. Furthermore, each data point results from averaging over
100 simulations, corresponding to 10 runs for each of 10 different
realizations of a given type of network of contacts with the
specified parameters. In the following simulations, $N=2001,
\kappa=2$ are kept invariant.

\begin{figure} \centering
\includegraphics[width=7.5cm]{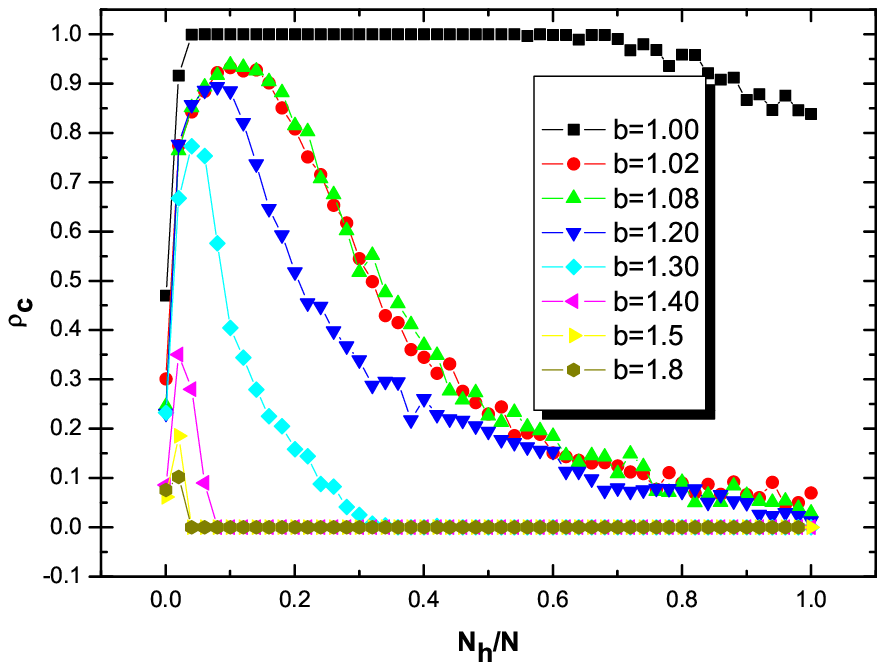}
\caption{\label{fig4}}
\end{figure}

\begin{figure} \centering
\includegraphics[width=7.5cm]{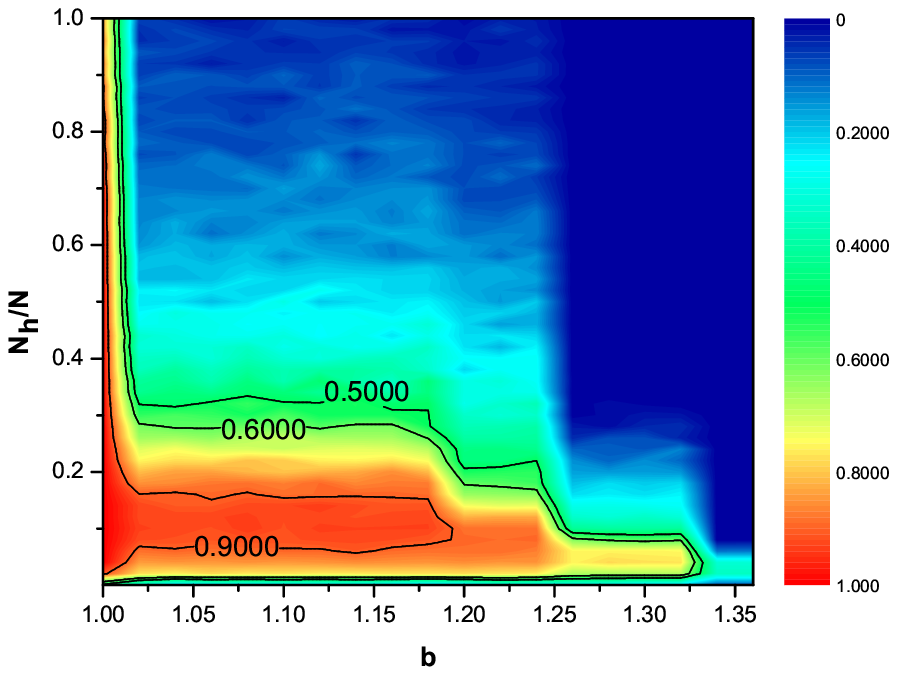}
\caption{\label{fig5} }
\end{figure}

To investigate the influence of network degree heterogeneity to
the evolution of cooperation, we fixed the value of $b$ when
increased fraction of hubs $N_h/N$ from $1/N$ to $1$.
Fig.~\ref{fig4} shows the equilibrium frequency of cooperators as
function of the fraction number of hubs $N_h/N$ by different
values of $b$ for comparison when $m=1000$. We find  the
$\rho_c$'s curve exhibits non-monotonous behavior with a peak at
some specific value of $N_h/N$. As aforementioned, the quantity
$N_h/N$ can be regarded as a measurement of heterogeneity. For
certain $b$, it has been illustrated in Fig.~\ref{fig4}, there is
a clear maximum frequency of cooperators around small $N_h/N$ near
the zero. In Fig.~\ref{fig5}, we report the result of frequency of
cooperation as function of the parameters space $(b,N_h/N)$. It
also shows that for fixed $b$, optimum cooperation level exits at
certain intermediate $N_h/N$. In addition, for certain
intermediate $N_h/N$, the density of cooperation is higher than
the cases when $N_h/N\to 0$ or $N_h/N \to 1$. As a result, it
reveals that optimal cooperation level occurs at intermediate
topological heterogeneity for some fixed $b$, that's neither the
most heterogeneous case nor the most homogeneous one favors the
cooperators best. Actually, under the most homogeneous case of
$N_h/N=1$, the cooperation behavior is almost extinct for large
values of $b$. While under the most heterogeneous case of
$N_h/N=1/N$, the cooperation level is also low. The reason for
this is that, in the most heterogeneous case, all shortcuts are
linked to a single node, making the network star-like with sparse
connections between the neighbor sites of the centered node. The
peripheral nodes' imitation strategies heavily depend upon the
centered individual's choice. Once the unique largest hub is
occupied by a defector, defection will spread over the entire
network. Then the system can hardly recover from the worst case of
all defectors. Therefore the most heterogeneous case with
$N_h/N=1$ is not the ideal case for arising cooperation. The most
homogeneous case with $N_h/N=1$ does not favor cooperation either.
In this situation, all nodes almost have the same number of
neighbors with some long range shortcuts making the average length
short. Therefore, the defection of some node is easily spread from
one node to another through the shortcuts. Within some
generations, the defectors prevail on the network. Hence the most
homogeneous case with $N_h/N=1$ does not promote the evolution of
cooperation.

\begin{figure} \centering
\includegraphics[width=7.5cm]{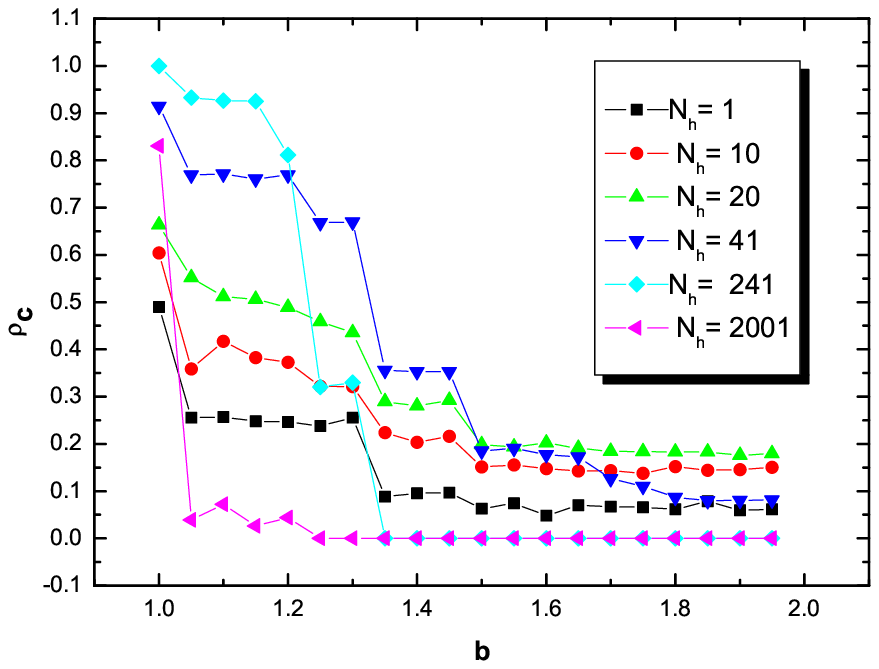}
\caption{\label{fig6}}
\end{figure}

Let's consider the case of intermediate heterogeneity where some
certain number of nodes are selected as hubs (density of
cooperators peaks at some specific value $N_h/N$ for fixed $b$).
These hubs are connected to each other by shortcuts or placed on
the same ring linked by the regular edges. This provides the
protection to each other in resisting the invasion of defectors.
Even if one hub is occupied by the defector, the cooperation level
will temporarily decrease due to the diffusion of defection from
the hub. However, other hubs occupied by cooperators could
reciprocally help each other and promote cooperation behaviors on
the network. Finally, cooperation is recovered to the normal
level. Accordingly, it in part result from the ability of
cooperators taking advantage of hubs. It is reported
in~\cite{Santos2005} that cooperation dominates in
Barab\'{a}si-Albert (BA) scale-free network due to the degree
heterogeneity where cooperators are more capable to occupy the
hubs in the network. Though the cooperation level on our HNW model
is not as high as it on the BA scale-free network (see
Fig.~\ref{fig6} and the figure in ref.~\cite{Santos2005} for
comparison), both of them benefit from the heterogeneity of the
network in the evolution of cooperation.

\begin{figure} \centering
\includegraphics[width=7.5cm]{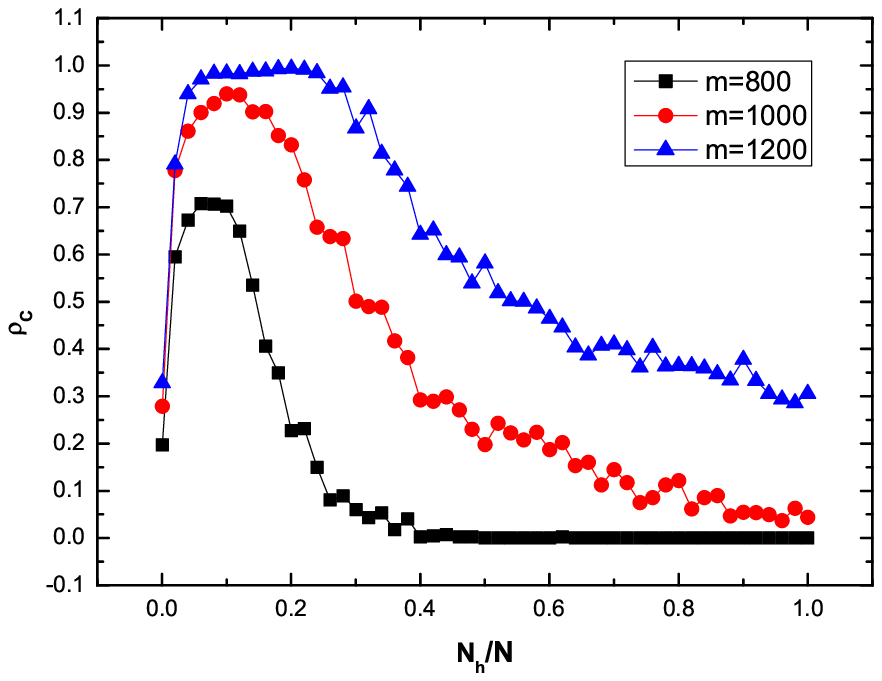}
\caption{\label{fig7}}
\end{figure}

The average degree of HNW network $\langle k
\rangle=2\kappa+2m/N$. It is believed that average degree affects
the evolution of cooperation. Fig.~\ref{fig7} plots the frequency
of cooperators versus $N_h/N$ corresponding to different shortcuts
number $m=800, 1000, 1200$ respectively with $b=1.1$. We find that
increasing average degree (larger $m$) promotes cooperation on the
network. Nevertheless, we should point out that with sufficient
large $m$ which makes the network nearly fully-connected, the
cooperation will be inhibited due to mean-field
behavior~\cite{wwx}. On the other hand, when $m$ changes, our
above conclusion on the effect of heterogeneity to cooperation is
still valid. We can see in Fig.~\ref{fig7} that, for different
fixed $m$ and $b$, $\rho_C$ still peaks at certain intermediate
value of $N_h/N$.

We further explore the situation where we adopt another
microscopic update rule different from Eq.~\ref{transp}. In stead
of updating the strategies according to the accumulated payoff of
individuals, we consider using the average payoff of individual
$x$, $\overline{P}_x=\frac{P_x}{k_x}$. Thus the new update rule
is:
\begin{equation}
W_{s_x \leftarrow s_y} = \frac{\overline{P}_y-\overline{P}_x}{b}
\end{equation}
By such update rule, we report the results in Fig.~\ref{fig8}
corresponding to $b=1.1$ and $b=1.3$. It is found that the
equilibrium frequency of cooperation is not sensitive to the
heterogeneity as $N_h/N$ increases from $1/N$ to $1$ (the curves
of average payoff are almost flat). It is because that normalizing
individual's payoff diminish the role of heterogeneity.
Consequently, the difference in individual's accumulated payoff
arising from the degree heterogeneity will be positive to
evolution of cooperation. We should be careful to come to a
conclusion with different update rules. Moreover, all these
phenomena and conclusions are also valid under the cases for
different population size $N$.

\section{Conclusion remarks and future work}
In conclusion, we have studied the effect of topological
heterogeneity to the evolution of cooperation behavior in
evolutionary Prisoner's Dilemma game. It is found that frequency
of cooperation peaks at some specific value of $N_h/N$, that is,
neither the most heterogeneous case nor the most homogeneous one
would favor cooperators. We found that the network degree
heterogeneity is one of the factors affecting the emergence of
cooperation. Besides, the average degree and the adopted update
rule also play a role to cooperation. Therefore, it is meaningful
and necessary to explore the underlying factors that do matter the
emergence of cooperation behaviors. In addition, the interplay
between game dynamics and network topologies is an interesting
topic for further investigations.

\section*{Acknowledgments} F.~F. would like to acknowledge Jie Ren for stimulating discussions.
The authors are partly supported by National Natural Science
Foundation of China under Grant Nos.10372002 and 60528007,
National 973 Program under Grant No.2002CB312200, National 863
Program under Grant No.2006AA04Z258 and 11-5 project under Grant
No.A2120061303.

\newpage

Fig.1 Illustration of the heterogeneous Newman-Watts small-world
model. A one-dimensional lattice with connections between all
vertex pairs separated by $\kappa$ or fewer lattice spacing, with
$\kappa=2$ in this case. The red dots are chosen as centers or
hubs and the added $m$ shortcuts (the dashed line) have at least
one endpoint belonging to these hubs.\\

Fig.2 The log-log plot of the variance of network degree sequence
versus fraction of hubs $N_h/N$ for $N=2001,m=1000$. Each data
point is averaged over 10 times different realization
corresponding to each $N_h/N$.\\

Fig.3 The plot of frequencies of cooperators as function of
temptation to defect $b$ corresponding to the different initial
fractions of cooperators $20\%,\,50\%,\,80\%$, respectively.
$N=2001, \kappa=2, m= 1000, N_h=41$.\\

Fig.4 Frequency of cooperators as function of the fraction of hubs
$N_h/N$, for different values of the temptation to defect $b$ as
shown in the legend. Note that when $b=1$, the game is not a
proper PD game.\\

Fig.5 The frequency of cooperators $\rho_c$ as function of the
parameter space $(b,N_h/N)$. The red region indicates the
cooperation of high level corresponding to the intermediate $N_h$
and $b$.\\

Fig.6 The frequencies of cooperators $\rho_c$ vs. the temptation
to defect $b$ with different levels of heterogeneity.\\

Fig.7 Plot of the fraction of cooperators vs. $N_h/N$ for
different values of $m$ with $b=1.1$.\\

\begin{figure} \centering
\includegraphics[width=7.5cm]{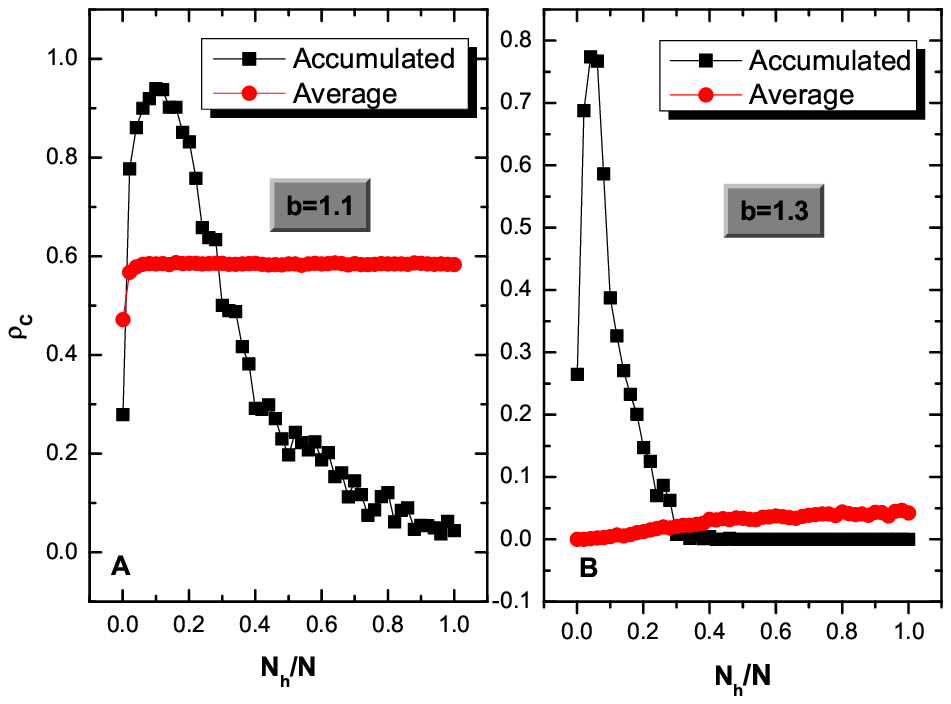}
\caption{\label{fig8}
 }
\end{figure}

Fig.8 Plot of the results with update rules corresponding to
average and accumulated payoff respectively. Fig.~\ref{fig8}-A
plots the case when $b=1.1$ and Fig.~\ref{fig8}-B with $b=1.3$\\

\newpage

\end{document}